\documentclass{elsart}
\usepackage{amsfonts}
\usepackage{graphicx}
\def\O{{\mathcal{O}}}
\def\g{{\mathfrak g}}
\def\t{{\mathfrak t}}
\def\X{{\mathfrak{X}}}
\def\R{{\mathbb R}}
\def\Xtan{\mathfrak{X}_{\rm tan}}
\def\Xt{X_{\rm tan}}
\def\Xi{X_{\rm inv}}
\def\tr{\mathop{\rm tr}}
\def\dbar{\overline{\nabla}}
\newcommand{\ddto}{{\left. \frac{d}{dt}\right |_{t=0}}}

\newtheorem{theorem}{Theorem}

\newtheorem{definition}{Definition}
\newtheorem{example}{Example}

\def\eqalign#1{\null\,\vcenter{\openup\jot \mathsurround=0pt \ialign{\strut
     \hfil$\displaystyle{##}$&$ \displaystyle{{}##}$\hfil \crcr#1\crcr}}\,}

\def\leqaligno#1{\displ@y \tabskip=\@centering \halign to\displaywidth{\hfil
$\@lign\displaystyle{##}$ \tabskip=0pt &$\@align\displaystyle{{}##}$\hfil
\tabskip=\@centering &\kern-\displaywidt\rlap{$\@lign##$}\tabskip=\displaywidth
\cr\cr #1\cr\cr}}

\begin{document}
\begin{frontmatter}
\title{Lie group foliations: \\
Dynamical systems and integrators}
\author{Robert I. McLachlan}
\address{IFS, Massey University, Palmerston North, New Zealand}
\author{Matthew Perlmutter}
\address{Departamento de Matematica, Instituto Superior Tecnico, Lisboa, Portugal}
\author{G. Reinout W. Quispel}
\address{Mathematics Department, La Trobe University, Bundoora, VIC 3086, Australia}

\begin{abstract}
Foliate systems are those
which preserve some (possibly 
singular) foliation of phase space, such as systems with integrals,
systems with continuous symmetries, and skew product systems. 
We study numerical integrators which also preserve the foliation.
The case in which the foliation is given by the orbits of an action of a Lie group
has a particularly nice structure, which we study in detail, giving conditions under
which all foliate vector fields can be written as the sum
of a vector field tangent to the orbits and a vector field invariant
under the group action. This allows the application of many techniques
of geometric integration, including splitting methods and Lie group integrators.
\end{abstract}
\end{frontmatter}

\section{Introduction}
In the early works of Feng Kang \cite{feng}, geometric integration was taken to be
the approximation of flows by elements of certain subgroups of $\mathit{Diff}(M)$ 
(the group of diffeomorphisms of the phase space $M$)---the groups of
symplectic, volume-preserving, or contact diffeomorphisms, for example.
This point of view was developed further in \cite{mc-qu} using the
Cartan classification of diffeomorphism groups and is continued here
by considering the so-called nonprimitive groups, those that leave
a foliation of $M$ invariant. (A preliminary announcement of
some of our results appears in \cite{an}.) That is, we ask: given a vector field whose flow
preserves a given foliation of $M$, how can we construct integrators
with the same property? Many other interesting (and difficult!) 
questions can be asked 
about this class of systems. How can the existence of an invariant
foliation be detected? What are the consequences for the dynamics
of the system? Even regarding
the construction of integrators, the class of all foliations seems to be too
large to admit a useful theory, and we are led (following the example
of Lie group integrators \cite{iserles,mu-za}) to consider foliations defined by
the action of a Lie group. We introduce these with an example.

\begin{example} \rm  Let $M = \R^2$ and consider the vector field
\label{ex:2d}
\begin{equation}
\label{eq:polar1}
\dot x = x y + x(1-x^2-y^2),\quad
\dot y = -x^2 + y(1-x^2-y^2) 
\end{equation}
In polar coordinates, this becomes
$$ \dot r = r(1-r^2),\quad \dot \theta = - r\cos\theta,$$
showing that the foliation into circles $r=$ const. is invariant
under the flow. (In fact, this foliation is {\em singular}, because the
leaf through the origin, a single point, has less than maximal dimension.)
A one-step integrator preserves this foliation (i.e., is `foliate'
 if the final value of $r$ is independent
of the initial value of $\theta$. Of course this is easy to
obtain in polar coordinates, but we shall see that
no standard integrator in cartesian coordinates is foliate.
The leaves of this foliation are the group orbits of the standard
action of $SO(2)$ on $\R^2$ (see Figure 1, top).
\end{example}

\begin{figure} \begin{center}
\vbox{\includegraphics[width=8cm,keepaspectratio,clip]{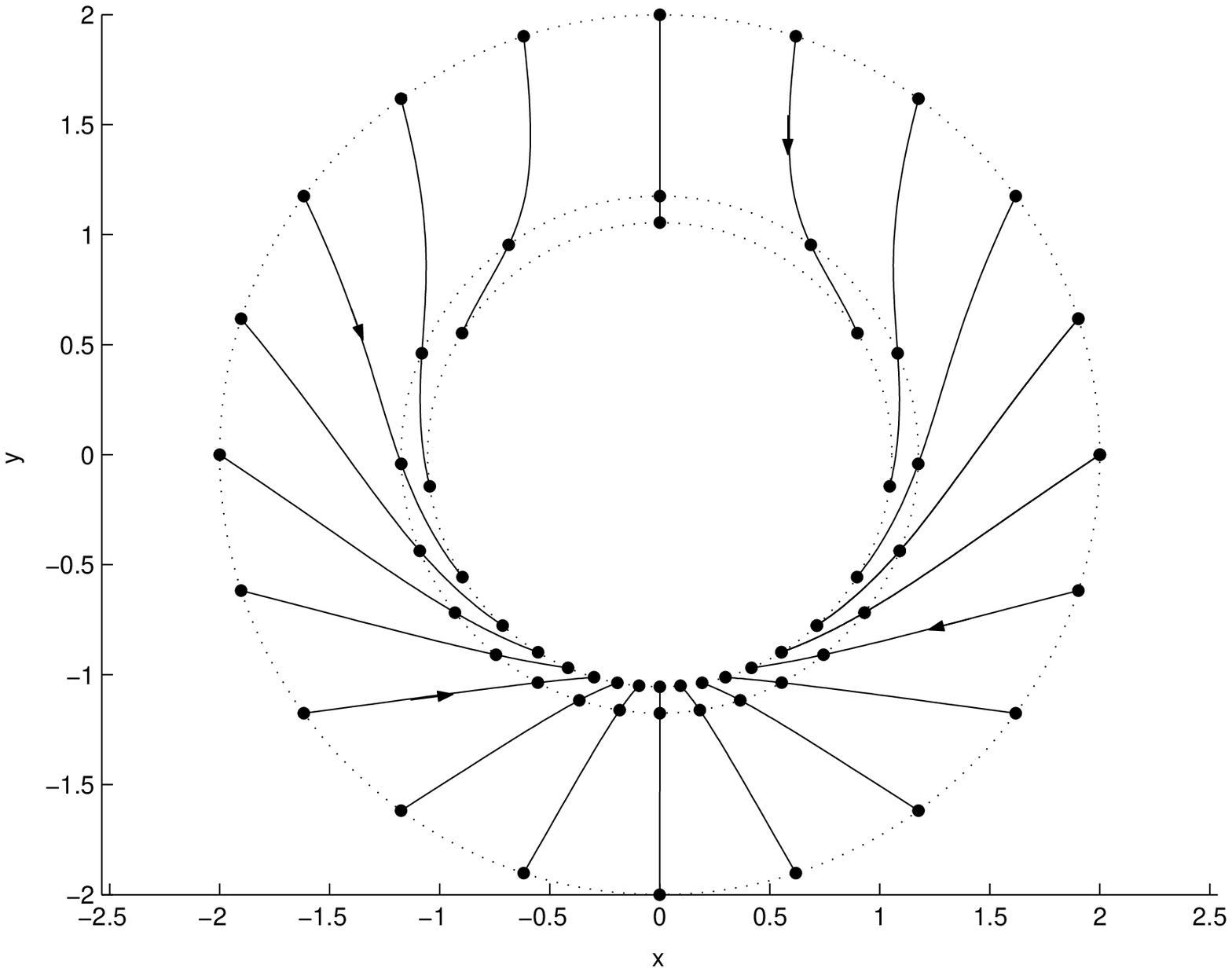} }
\bigskip
\vbox{\includegraphics[width=8cm,keepaspectratio,clip]{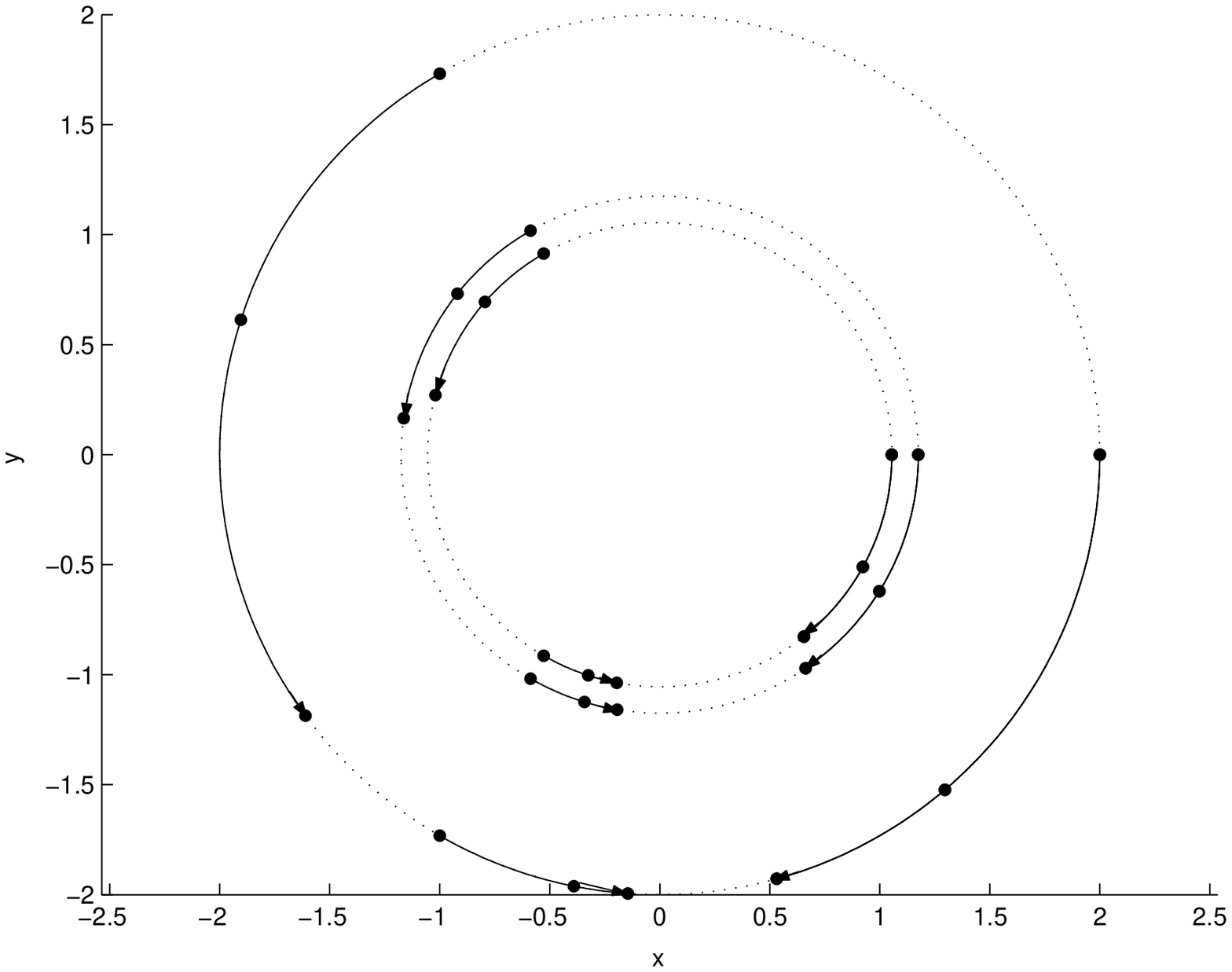} }
\bigskip
\vbox{\includegraphics[width=8cm,keepaspectratio,clip]{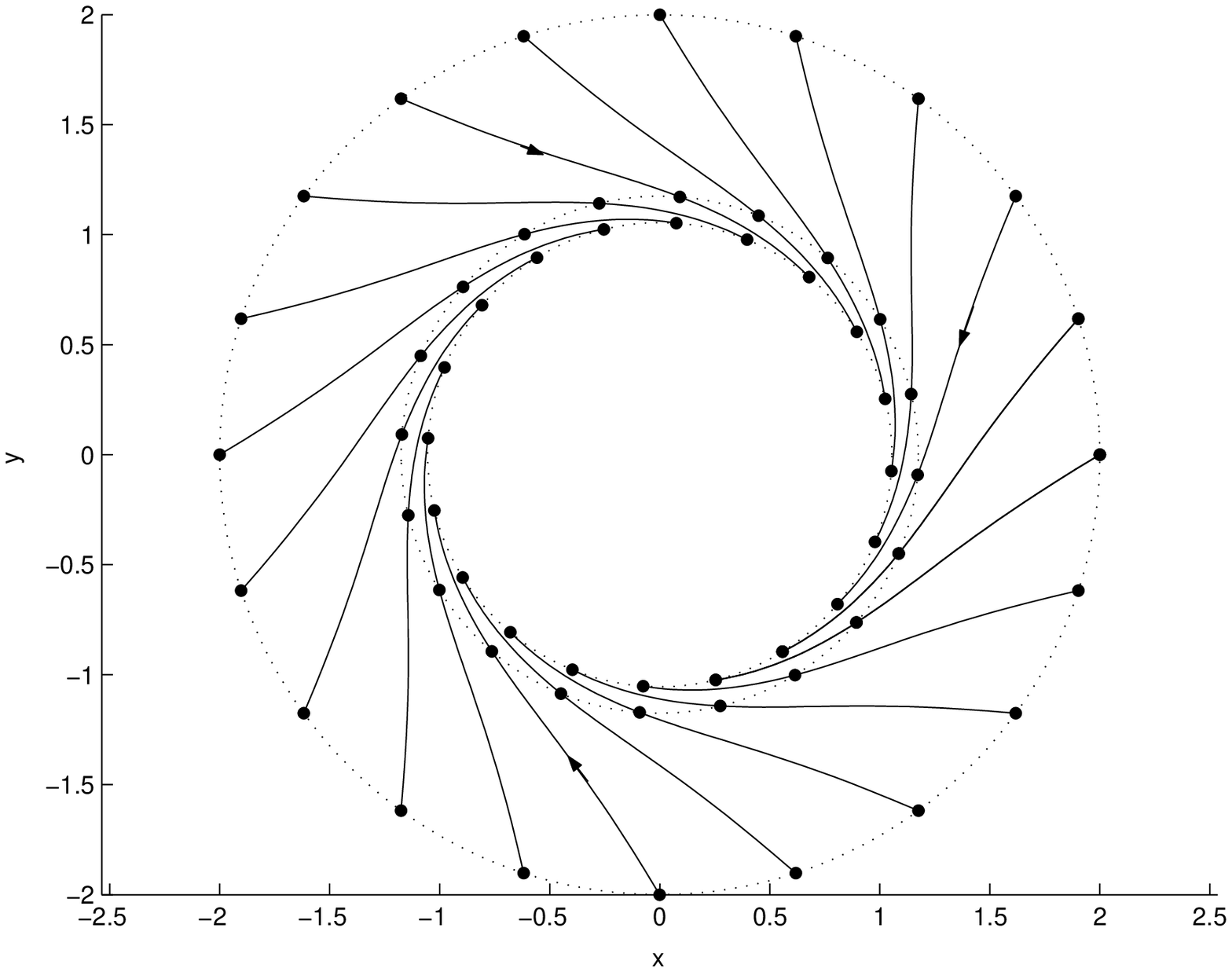} }
\caption{Three foliate vector fields.
Top: a general foliate vector field,
$\dot r = r(1-r^2)$, $\dot\theta = -r\cos\theta$ (Eq. (\ref{eq:polar1})). Middle: a system
with an integral, $\dot r = 0$, $\dot\theta = -r\cos\theta$.
Bottom: a system with a continuous symmetry, $\dot r = r(1-r^2)$,
$\dot\theta = -(1+r^2/5)$. All three flows map circles to circles. The
dots mark times 0, 0.5, and 1.}
\label{foliate.eps}
\end{center} \end{figure}

We show in section \ref{sec:liegroup} that in many cases, vector fields whose flows
preserve such a foliation have a particularly simple structure, namely,
they can be written as a sum of two vector fields, one tangent to the
leaves of the foliation (being
the orbits of a group action), and one invariant under the group action. 
This representation allows the construction of various types of
foliation-preserving integrators.

\begin{definition} \cite{molino,sussmann}
Let $M$ be a manifold of dimension $m$.
A {\em singular foliation} $F$ of $M$ is a partition of $M$ into
connected immersed submanifolds (the ``leaves''), such that the vector
fields on $M$ tangent to the leaves are transitive on each leaf
\cite{sussmann}. $F$ is {\em regular} if each leaf
has the same dimension. $F$ has {\em codimension $q$} if the maximum dimension
of the leaves of $F$ is $m-q$.
A diffeomorphism of $M$ is {\em foliate} with respect to $F$
if it leaves the foliation invariant, i.e., if it maps leaves to leaves.  
A vector field on $M$ is {\em foliate} if its flow is foliate. 
The space of smooth vector fields foliate with respect to $F$ is denoted $\X_F$.
A (one-step) integrator is {\em foliate with respect to $F$ }
if the diffeomorphisms of $M$ corresponding 
to each time step are foliate.
The space of smooth vector fields tangent to the leaves of $F$ is 
denoted $\Xtan$. 
The {\em space of leaves} (denoted $M/F$) is obtained by identifying
the points in each leaf together with the quotient topology.
\end{definition}

\begin{theorem} \cite{molino}
\label{prop:general}
$\X_F$ and $\Xtan$ form Lie algebras. $\Xtan$ is
an ideal in $\X_F$. A vector field $X$ is foliate with respect to $F$
if and only if $[X,Y]\in \Xtan$ for all $Y\in \Xtan$. 
\end{theorem}


\begin{theorem} \cite{molino}
\label{prop:simple}
Let $M$ and $N$ be manifolds  of dimension $m$
and $n$, respectively. Let $I:M\to N$ be a smooth surjection. (If
$I$ is not onto, we replace $N$ by $I(M)$.) Then
$I$ defines a foliation $F$ whose leaves are given by the connected
components of $I^{-1}(y)$ for each $y\in I(M)$. If $I$ is a submersion,
i.e. if $dI$ has constant rank $n$, then $F$ is a regular foliation
of codimension $n$. In this case the space of leaves $M/F$ is diffeomorphic
to $N$. 
\end{theorem}

Such a foliation is called {\em simple}. Given a vector field, one can
search for simple foliations it preserves by looking for functions $I$
such that $\dot I = f(I)$. These functions $I$ might be conserved
quantities in special cases. For example, the conserved momentum 
function $J$ of a Hamiltonian system can evolve under
$\dot J = c J$ under the addition of linear friction. Then the
perturbed system preserves the simple foliation defined by $J$ \cite{mc-pe}.

For a regular simple foliation, a foliate vector field $X$ drops
to a vector field on the space of leaves $M/F$, which (following
the tradition in geometric mechanics, \cite{ma-ra}) we call
the {\em reduced} system. If $\xi:[0,T]\to M/F$ is an orbit
of the reduced system, 
the problem of finding an orbit of $X$ is called the {\em
reconstruction problem}. If all leaves are diffeomorphic, then reconstruction
involves integrating a nonautonomous vector field on a fixed leaf.
(In Example \ref{ex:2d}, the reduced system is $\dot r = r(1-r^2)$;
for any solution $r(t)$ of this equation, the reconstruction system
is $\dot \theta = - r(t)\cos\theta$, a nonautonomous ODE on the leaf
$S^1$.) 

In the application to integrators, we do not usually want to construct
the reduced system explicitly since the original phase space $M$ is usually
linear and easier to work in. We would like integrators on $M$ which
preserve the foliation automatically.

\begin{example} \rm A system with $k$ first integrals $I:M\to\R^k$ is foliate
with respect to the level sets of the functions $I$. Each leaf is in fact
fixed by the flow. For this reason we choose the symbol $I$ in Theorem
\ref{prop:simple} to suggest that simple foliate systems generalize systems
with first integrals. (See Figure 1, middle.)
\end{example}

\begin{example} \rm A system with a symmetry is foliate with respect to
the orbits of the symmetry. That is, let $X$ admit the Lie group action
$\lambda:G\times M\to M$ as a symmetry, so that its flow $\phi_t$ 
is $G$-equivariant. 
Then $\lambda(g,\phi_t(x)) = \phi_t(\lambda(g,x))$, i.e., the foliation
with leaves given by the group 
orbits $\{\lambda(g,x):g\in G\}$ is invariant.
In this case the reconstruction problem on $G$ is easier to solve
than in the general case, because it is $G$-invariant. (See Figure 1, bottom.)
\end{example}

\begin{example} \rm
\label{ex:lorenz}
The Lorenz system is given by
$$\eqalign{\dot x &= \sigma y - \sigma x, \cr
\dot y &= -y - x z - r x, \cr
\dot z &= x y - b x. \cr
}$$
If $b=2\sigma$, the system is foliate with leaves $x^2 - 2\sigma z=$ const., for
$$ \frac{d}{dt}(x^2-2\sigma z) = -2\sigma(x^2 - 2\sigma z).$$
(A foliate integrator can be constructed as follows. We split into
$$\begin{array}{llll}
X_1\colon& \dot x = \sigma y,& \dot y = -x z - r x,& \dot z = x y, \\
X_2\colon& \dot x = -\sigma x,& \dot y = -y,& \dot z = -2\sigma x.\\
\end{array}
$$
$X_1$ is tangent to the foliation and may be integrated using the midpoint
rule, which preserves the quadratic function $x^2-2\sigma z$;
$X_2$ is foliate but linear, and can be solved exactly. Composition
then yields a foliate integrator.)
\end{example}

\begin{example} \rm
\label{ex:product}
A special case of the foliations defined by submersions is given
by $M=N\times L$ and $I$ is projection onto $N$. Each leaf is 
then diffeomorphic to $L$. (This is locally true of any simple foliation 
in a neighborhood in which the leaves have constant dimension.)
In coordinates
$x$ on $N$ and $y$ on $L$, any foliate vector field can be written
in coordinates as
$$ \eqalign{ \dot x &= f(x) \cr
             \dot y &= g(x,y). \cr
}$$
and any tangent vector field as
$$ \eqalign{ \dot x &= 0 \cr
             \dot y &= g(x,y). \cr
}$$
These foliate vector fields are also known as {\em skew product systems},
introduced by Anzai \cite{anzai}, and studied today in ergodic theory 
\cite{bu-wi} and complex dynamics \cite{jonsson}.
A special case of skew product systems is given by $N=\R^n$ and
$L=\R^{m-n}$, the standard foliation on $\R^m$.
\end{example}

\begin{example} \rm The extension of a nonautonomous vector field on $M$
to an autonomous vector field on $M\times\R$ preserves the foliation 
defined by $t=$ const. Most integrators
are foliate and indeed, solve the reduced system $\dot t = 1$ exactly.
\end{example}

In a foliate system, one can obtain some information about part of the
system (namely, the current leaf) for all time without even knowing the
full initial condition. Surely this puts strong dynamical constraints
on the whole system. Nevertheless, the only properties we can point to
depend on the reduced or reconstruction systems being simple in some
way: having small dimension or simple (e.g. linear) dynamics.

\begin{example} \rm Burns and Wilkinson \cite{bu-wi} study skew product
systems of the form $(x,g)\mapsto (f(x),\phi(x)g)$ on $M\times G$, where
$G$ is a compact Lie group, $\phi:M\to G$, and $f$ is measure preserving.  Here the
foliation is simple as in Example \ref{ex:product} and the systems
have extra structure: the reduced system $x\mapsto f(x)$ is measure
preserving, the reconstruction system $g\mapsto \phi(x)g$ is $G$-equivariant,
and the whole system is measure preserving (with respect to the product
of the measure on $M$ with Haar measure on $G$). This class is a
special case of the class of measure preserving systems with a symmetry. The
extra structure allows them to prove that the ergodic systems are
open and dense in this class.
\end{example}

\begin{example} \rm Consider the system on $\R^3$ with a codimension 1
foliation,
$$ \dot x = f(x),\ \dot y = g(x,y,z),\ \dot z = h(x,y,z).$$
The only possible $\omega$-limit set of the reduced system $\dot x = f(x)$
is a point, which suggests that in this case the $\omega$-limit set of the whole system is either
a point or a circle (periodic orbit). 
(If $f(x) > 0$, however, so that the reduced system has
no $\omega$-limit set, chaos is possible in the full system.)
Similarly, for the system on $\R^3$ with a codimension 2 foliation,
$$ \dot x = f(x,y),\ \dot y = g(x,y),\ \dot z = h(x,y,z).$$
the possible $\omega$-limit sets of the reduced system in $(x,y)$ are a points
and circles, which suggests that in these cases the $\omega$-limit set 
of the full system is a point,a circle, or a 2-torus.
In both examples
the existence of the foliation influences the possible dynamics.
\end{example}

\section{Integrators for simple foliations}
In this section we consider whether standard integrators can be foliate,
and adapt geometric integrators for systems with integrals to the foliate
case.

The next Theorem generalizes the fact that Runge-Kutta methods preserve
arbitrary linear integrals.

\begin{theorem}
Let $M=\R^m$ and let $F$ be a linear foliation, i.e., a simple foliation
defined by the linear function $I:\R^m\to \R^k$. Then any 
Runge-Kutta method is foliate.
\end{theorem}
\begin{pf}
Runge-Kutta methods are linearly covariant \cite{mc-qu-tu}, hence we can
apply a linear change of variables to bring $I$ into the form
$I(x_1,\dots,x_m) = (x_1,\dots,x_k)$. Runge-Kutta methods are closed
under restriction to the closed subsystem with coordinates $x_1,\dots,x_k$
\cite{bo-sc}, hence the final values of $x_1,\dots,x_k$ depend only
on their initial values, i.e., the method is foliate.
\end{pf}

We now give two methods for constructing foliate integrators, each based on
a popular integral-preserving method.

\noindent{\bf Projection methods.} Let the leaves of the foliation
be the level sets of the functions $I:M\to\R^k$ and write the reduced
system as $\dot I = h(I)$. Let $x_n$ be the initial condition.
Then the following algorithm provides a foliate integrator.
\begin{list}{}{\leftmargin=1.5cm \labelwidth=1.5cm}
\item[Step 1.] Calculate $I_n = I(x_n)$ and integrate the reduced
system for one time step, giving $I_{n+1}$, which depends only on
$I_n$.
\item[Step 2.] Integrate the full vector field for one time step,
giving $\tilde x$.
\item[Step 3.] Calculate $x_{n+1}$ by projecting $\tilde x$ onto
the desired leaf $I^{-1}(I_{n+1})$, e.g. by orthogonal
projection.
\end{list}

\noindent{\bf Discrete gradient methods.} \cite{mc-qu-ro}
For simplicity we present the method for a single function $I$,
i.e., for a codimension 1 foliation, and use a Euclidean metric. 
It is necessary to first
split the vector field into its components tangent and orthogonal
to the leaf. Since $\dot I = X.\nabla I = h(I)$, the component of $X$ 
in the direction $\nabla I$ is $(h(I)/|\nabla I|^2)\nabla I$, while
the component of $X$ tangent to the leaf can be written in the form
$A(x)\nabla I$ for some antisymmetric matrix $A(x)$ \cite{mc-qu-ro}.
That is, the full system can be written
$$\dot x = \left(A(x) + \frac{h(I)}{|\nabla I|^2}\right)\nabla I.$$
Extend $A$ and $h$ to functions $\bar A$ and $\bar h$ of pairs of points
satisfying $\bar A$ antisymmetric,
$\bar A(x,x)=A(x)$, $\bar h(I,I) = h(I)$, and let $\dbar$ be
a discrete gradient. Then
the discrete-gradient discretization
$$ \frac{x'-x}{\tau} = \left(\bar A(x,x') + \frac{\bar h(I(x),I(x'))}
{|\dbar I(x,x')|^2}\right)\dbar I(x,x') $$
obeys
$$ I(x') - I(x) = \tau \bar h(I(x),I(x')),$$
so the integrator is foliate. Here $\tau$ is the time step.

One popular integral-preserving method does not generalize to foliations.
The symplectic Runge-Kutta methods preserve arbitrary quadratic first
integrals (and arbitrary linear symmetries) but do not preserve arbitrary
quadratic foliations. To see this, consider the midpoint rule applied
to a system preserving the foliation $r=$ const. as in Example 1. Let
the method be $x\mapsto x'$ and write $\bar x = (x+x')/2$. The method
is $x' = x + \tau X(\bar x)$. We have
$$ {r'}^2-r^2 = (x'+x)^t (x'-x) = 2\tau \bar x^t X(\bar x).$$
Since the vector field $X$ is foliate, 
$$\frac{d}{dt}\frac{1}{2}r^2 = x^t \dot x = x^t X(x) $$
is a function of $r$ only. Therefore, 
${r'}^2 - r^2$ is a function of $|\bar x|^2$ only. However, $|\bar x|^2$
is not in general a function of $r$.
For example, consider the system
\begin{equation}
\label{eq:polar2}
\dot x = -y^2 + x,\quad \dot y = x y + y
\end{equation}
for which $\dot r = r$, $\dot\theta = r\sin\theta$. Applying the midpoint
rule and expanding in a Taylor series gives
$$ {r'}^2 = r^2(1+2\tau + 2\tau^2 + \tau^3(3-y^2)/2)+\O(\tau^4) $$
which is not a function of $r$ only.

To obtain foliate integrators which are `intrinsic' in the sense
that they do not involve constructing the reduced system explicitly,
we have to consider special cases. For example, linear foliations can
be preserved in a sense automatically, and for systems with linear symmetries,
the foliation defined by the orbits of the symmetry is preserved by
symmetry-preserving integrators.
We therefore consider in the next section foliations defined by the
orbits of a Lie group action. Not only is the foliation simple, but
the foliate vector fields themselves have a nice structure which allows
the construction of foliate integrators. 

\section{Lie group foliations}
\label{sec:liegroup}

Let $G$ be a Lie group and $\lambda\colon G\times M\to M$ be an action of
$G$ on $M$. (We write the
action as 
$\lambda(g, x) = \lambda_g x = g\cdot x = gx$, as needed.)
This group action generates a (possibly singular) foliation
whose leaves are the group orbits $\lambda(G,x)$. A vector field
preserving this foliation is said to be {\em $G$-foliate}.
Let $\g$ be the Lie
algebra of $G$ and let $\g_M$ be the distribution tangent to the leaves, with each
$\xi\in\g$ associated to a vector field $\xi_M$ on $M$, i.e.,
$$ \g_M(x) = T_x(\lambda(G,x)) = \{\xi_M(x)\colon  \xi\in \g\}.$$ 
(Recall that $p$-dimensional distribution on $M$ is given by associating, to each
point $x$ in $M$, a $p$-dimensional subspace of the tangent space $T_x M$ 
\cite{molino}.)

As with any foliation, any vector field tangent to the leaves is foliate. 
However, in the case of a Lie group foliation we have another natural 
class of foliate vector fields, namely those invariant under the action. 
For, from Theorem 1, $X$ is foliate iff $[X, \g_M]\subset \g_M$; 
but if $X$ is $G$-invariant, then $[X , \g_M] = 0$. Thus, the vector field 
$\Xt+\Xi$ is $G$-foliate, where $\Xt$ is tangent to the leaves and
$\Xi$ is $G$-invariant. To characterize the $G$-foliate vector fields we 
now want to turn this around and ask: When can all foliate vector fields 
be decomposed in this way? We shall see that this decomposition makes it 
easy to construct foliate vector fields explicitly and also to construct 
foliate integrators. 

Our results can be summarized as follows. 
Let X be a foliate vector field. First, if $G$ acts by isometries, 
then $X = \Xt+\Xi$ where $\Xi$ is not only $G$-invariant but also 
perpendicular to the leaves (Theorem 5). Second, if $\g_M$ has a
$G$-invariant complement $H$, then $X = \Xt+\Xi$ with $\Xi\in H$ 
(Theorem 6). This arises, for example, if the action is free and proper 
(Theorem 7). The decomposition is also true locally in a neighborhood of an 
orbit which admits a certain slice (Theorem 8). Finally, we give a 
counterexample to show that not all actions admit such a decomposition 
(Example 9).

Recall that a proper action is one for which $G\times M \to M \times
M$, $(g,x)\mapsto (g.x,x)$ is proper, i.e., the inverse image
of any compact set is compact. An action of a compact group
must be proper. A proper action must have compact isotropy
groups at all points of $M$.

\begin{theorem}
\label{theorem:proper}
\cite{du-ko} Let $G$ be a Lie group acting properly on the manifold $M$. 
Then $M$ has a $G$-invariant Riemannian metric.
\end{theorem}

\begin{theorem}
\label{theorem:main}
Let $(M,\left\langle,\right\rangle)$ be a Riemannian manifold. Suppose we
have a smooth group action $\lambda\colon G\times M\rightarrow M$ acting by
isometries on the metric $\left\langle,\right\rangle$. Consider
the foliation by the group action where the leaf through $x\in M$
is simply the orbit $\lambda(G,x)$. 
Then the vector field $X$ is $G$-foliate if and only if
the unique metric decomposition of $X$ into tangential 
and perpendicular components
\begin{equation}
X=X^\|+X^\perp
\end{equation}
satisfies, for all $g\in G$, $\phi_{g}^{\ast}X^{\perp}=X^{\perp}$.
In other words,  $X$ is $G$-foliate if and only if
its component perpendicular to the leaves is $G$-invariant.
\end{theorem}
\begin{pf}
Recall that $\Xtan$ is the space of smooth vector fields tangent 
to the foliation. In this case we can write
\begin{equation}
\Xtan(x)=\{\xi_{M}(x)\colon\xi\in\mathfrak{g}\}.
\end{equation}
From Theorem \ref{prop:general} we have that
$X$ is foliate if and only if for all $Y\in \Xtan$,
\begin{equation}
[X,Y] = [X^{\|}+X^{\perp},Y]\in\Xtan.
\end{equation}
However, $X^{\|}\in\Xtan$ by definition, so
$X$ is foliate if and only if for all $Y\in \Xtan$,
\begin{equation}
[X^{\perp},Y]\in \Xtan.
\end{equation}
Now, suppose $X$ is foliate. The condition above requires, in particular,
that for all $\xi\in\mathfrak{g}$,
$[\xi_M,X^{\perp}]\in\Xtan$.
Now, by definition of $X^{\perp}$ we have that $\left\langle X^{\perp},
\xi_{M}\right\rangle=0$. Fix $x\in M$ and $\eta\in \mathfrak{g}$. Let
$\phi_{t}^{\eta}$ denote the flow of $\eta_{M}$ which is given by
$\phi_{t}^{\eta}(x)=\lambda(\exp t\eta,x)$. We then have
\begin{equation}
\left\langle X^{\perp}(\phi_t^\eta(x),\xi_{M}(\phi_t^\eta(x))
\right\rangle=0.
\end{equation}
Since $\phi_{t}^{\eta}$ leaves the metric invariant for all $t$,
we have 
\begin{equation}
\left\langle T\phi^{\eta}_{-t}X^{\perp}(\phi^{\eta}_{t}(x)),
T\phi^{\eta}_{-t}\xi_{M}(\phi^{\eta}_t(x))\right\rangle=0
\end{equation}
so that
\begin{equation}
\label{eq:Lie}
\eqalign{0 &=
\ddto \left\langle T\phi^{\eta}_{-t}X^{\perp}(\phi^{\eta}_{t}(x)),
T\phi^{\eta}_{-t}\xi_{M}(\phi^{\eta}_t(x))\right\rangle \cr 
&=\left\langle L_{\eta_{M}}X^{\perp}(x),\xi_{M}(x)\right\rangle+
\left\langle X^{\perp}(x),L_{\eta_{M}}\xi_{M}(x)\right\rangle.
}
\end{equation}
Since the action is on the left we have
$L_{\eta_{M}}\xi_{M}=[\eta_{M},\xi_{M}]=-[\eta,\xi]_{M}$, so that
the second term in Equation (\ref{eq:Lie}) is $0$. It follows that
\begin{equation}
\left\langle L_{\eta_{M}}X^{\perp},\xi_{M}\right\rangle=0
\end{equation}
which shows that $L_{\eta_{M}}X^{\perp}$ is perpendicular to the
leaves. However, by the assumption that $X$ is foliate we demand that
$[X^{\perp},\eta_{M}]=L_{\eta_M}X^\perp\in\Xtan$ so that we must
have $L_{\eta_{M}}X^{\perp}=0$, proving $G$-invariance of $X^{\perp}$.
Conversely, suppose $X^{\perp}$ is $G$-invariant. It suffices to show
that $X^{\perp}$ is foliate with respect to the $G$ orbits. But this
is immediate since the flow of $X^{\perp}$ is $G$-equivariant
and therefore, if $\phi_{t}$ denotes the flow of  $X^{\perp}$,
$\phi_{t}\lambda(g,x)=\lambda(g,\phi_{t}(x))$ which shows that $\phi_{t}$
maps the leaf through $x$ to the leaf through $\phi_{t}(x)$.
\end{pf}

\begin{theorem}
\label{thm:comp}
 Let $\lambda\colon G\to M$ be a smooth action of the Lie group $G$ on the
manifold $M$. Suppose the tangent bundle $TM$ admits a $G$-invariant splitting
$TM = \g_M + H$. Let the corresponding decomposition of a vector field $X$ be
$X = Y + Z$. Then $X$ is $G$-foliate iff $Z$ is $G$-invariant.
\end{theorem}
\begin{pf}
The `if' part is immediate. Conversely, let $X$ be a $G$-foliate vector
field. We have to show that $Z$ is $G$-invariant, i.e., that $[\xi_M, Z] = 0$ for all $\xi\in\g$.
Since $X$ is foliate, $[\xi_M,X] = [\xi_M, Y + Z]\in \g_M$; but $[\xi_M, Y ]\in \g_M$ 
since $Y$ is 
tangent to the leaves by definition. Therefore $[\xi_M, Z] \in \g_M$.
Fix $x \in M$ and $g \in G$. Since $H$ is $G$-invariant, we have 
$$T_{g\cdot x}\lambda_{g^{-1}}(Z(g\cdot x))\in H_x.$$
Letting $g = \exp(t\xi_M)$ and differentiating with respect to $t$ at $t = 0$ gives
$$\eqalign{
[\xi_M,Z](x) &= L_{\xi_M} Z(x) \cr
&= \ddto \lambda^*_{\exp(t\xi)} Z(x) \cr
&= \ddto T_{\exp(t\xi)\cdot x} \lambda_{\exp(-t\xi)\cdot x}
(Z(\exp(t\xi)\cdot x))\in H_x.}$$
Since the two subspaces $\g_M(x)$ and $H_x$ are complementary, we must have
$[\xi_M, Z] = 0$, that is, $Z$ is $G$-invariant.
\end{pf}

The next theorem follows either from Theorem 5 or from Theorem 6. However,      
we also give a direct proof which constructs the decomposition explicitly
relative to the choice of a connection. Recall that a free action is one in which
all isotropy groups are trivial, $G_x = \{e\}$ for all $x\in M$.

\def\A{{\mathcal A}}
\def\hor{\mathop{\rm hor}}
\def\M#1{\left( #1 \right)_M }

\begin{theorem}
\label{thm:free}
Let $\lambda\colon G\times M \to M$ be a smooth, free, and proper action of the
Lie group $G$ on the manifold $M$. Let $X$ be a $G$-foliate vector field. Then there
exists a vector field $\Xt\in\Xtan$ and a $G$-invariant vector field $\Xi$ such that
$X = \Xt + \Xi$.
\end{theorem}
\begin{pf}
Because the action is free and proper, $M \to M/G$ is a principal
$G$-bundle. There is a 1--1 correspondence between the space of connections
and $G$-invariant complements to $\g_M$. Fix a connection $\A$ and let $H_x = 
\ker \A_x$.
Furthermore, there is a unique horizontal lift map $\hor\colon  T(M/G) \to TM$. It
then follows that $\hor \pi_*X$ is a $G$-invariant vector field taking values in $H$.
Thus $\Xi = \hor \pi_*X$ and $\Xt = X-\Xi$ provides the desired decomposition
of $X$. $\Xt$ is tangent to the foliation because $\A(\Xt) = \A(X - \hor \pi_*X) =
\A(X) - \A(\hor \pi_*X) = \A(X)$.
\end{pf}

\begin{definition} \cite{du-ko}
A slice at $x_0\in M$ of a smooth Lie group action $\lambda\colon G\times M \to
M$ is a submanifold $S$ of $M$ through $x_0$  such that    
\begin{itemize}
\item[{(i)}]
$T_{x_0} M = \g_M(x_0 ) \oplus T_{x_0} S$ and for all $x\in S$, $T_x M = \g_M(x) + T_x S$;
\item[{(ii)}]     
$S$ is $G_{x_0}$-invariant; and
\item[{(iii)}]
if $x \in S$, $g \in G$, and $\lambda(g, x) \in S$, then $g \in G_{x_0}$.
\end{itemize}
\end{definition}

Note that $G \cdot S$ is an open, $G$-invariant neighborhood of the orbit $G \cdot x_0$,
and every orbit in this neighborhood intersects $S$ in a unique $G_{x_0}$-orbit of $S$.
Every proper action admits a slice at every point. The dimension of an orbit
which intersects $S$ must be greater than or equal to the dimension of the orbit
through $x$.

Unfortunately, it appears that existence of a slice is not sufficient to guarantee
the existence of a $G$-invariant splitting of $T(G\cdot S)$. The following theorem
requires two extra assumptions. However, the first is a requirement on the
group itself (not on the action), while the second only concerns the action
of $G_{x_0}$ on $S$, reducing the dimensionality. One important special case (to be
illustrated in Example 10) is when $G_x = G_{x_0}$ for all $x$ in $S$, in which the second
assumption is automatically satisfied with $H = TS$. 
The following theorem   
gives sufficient conditions which, together with Theorem 6, guarantee that on
$G\cdot S$ every $G$-foliate vector field can be decomposed into tangent and invariant
components.

\begin{theorem}
\label{thm:slice}
Let $\lambda\colon G\times M \to M$ be a smooth action of the Lie group $G$ on the
manifold $M$ which admits a slice $S$ through the point ${x_0} \in M$. Furthermore,
assume
\begin{itemize}
\item[(i)] $\g$ admits a $G_{x_0}$-invariant splitting $\g = \g_{x_0} \oplus \t$ 
where the group $G_{x_0}$ acts on
$\g$ by adjoint action.
\item[(ii)] There exists a $G_{x_0}$-invariant splitting $TS = \M{\g_{x_0}} \oplus H$ of the tangent
bundle of $S$.
\end{itemize}
Let $U$ be the open $G$-invariant neighborhood $G\cdot S$ of the orbit through $ {x_0}$. Then 
the tangent bundle $TU$ admits a $G$-invariant splitting
\begin{equation}
\label{eq:slice}
T_u M = g  \cdot H_s \oplus \left[\M{g \cdot \t} (u) \oplus \M{\g_{g\cdot {x_0}}}
(u)\right],
\end{equation}
where $u \in U$ is given by $u = gs$, the first distribution on the right hand side is
transverse to the group orbits and the other two form a decomposition of the
group directions into directions transverse and tangent to the slice respectively.
\end{theorem}
\begin{pf}
We first check that each distribution in Eq. (\ref{eq:slice}) is $G$-invariant.   
Since the action satisfies $G_{g\cdot{x_0}} = g G_{x_0} g^{-1}$, 
we have $ \g_{g\cdot x_0} = \mathop{\rm Ad}_g \g_{x_0} =: g \cdot \g_{x_0}.$
Therefore
$$T \lambda_g\colon \M{\g_{x_0}} (s) \simeq \M{g\cdot  \g_{x_0}} (gs).$$
Next, since $g \cdot\M{\xi} =  \M{g \cdot \xi}$ holds for any group action, we have
$$T \lambda_g\colon \M{\t}(s) \simeq \M{g\cdot  \t} (gs).$$
It follows that $\M{g \cdot  \t} (gs) \cap \M{g \cdot  \g_{x_0}} (gs) = 0$. 
Next, the first distribution is   
$G$-invariant by construction, and since $T \lambda_g\colon T_s(G \cdot s) \simeq 
T_{gs}(G \cdot gs)$, we must
have $g \cdot H_s \cap T_{gs}(G \cdot gs) = 0$.

Finally, we must check that the splitting is well defined. Suppose $u = gs =
g_1 s_1$. Then without loss of generality we can take $s_1 = s$ and $g_1 = gh$ with   
$h \in G_{x_0}$, since by the slice property, $G_s \subset G_{x_0}$. 
(Each orbit in $U$ intersects $S$
in a unique $G_{x_0}$-orbit.) The splitting is then well-defined if $T \lambda_h$ maps each of
the three distributions isomorphically into themselves. We clearly have that
$T \lambda_h\colon  T_s S \simeq T_{hs}S$ and that 
$T \lambda_h\colon  T_s(G_{x_0} \cdot s) \simeq T_{hs}(G_{x_0} \cdot s)$, since 
$h \in G_{x_0}$ . So
what we need is that $T \lambda_h$ maps $H_s$ into $H_{hs}$, which is true 
by assumption (ii),
and that $h  \cdot \t = \t$ for all $h \in G_{x_0}$ , which is true by assumption (i).
\end{pf}

We do not have necessary and sufficient conditions for the action to be such
that the tangent and invariant vector fields span all foliate vector fields. How
ever,
the following example shows that this is not true for all group actions.

\def\dx{\frac{\partial}{\partial x}}
\def\dy{\frac{\partial}{\partial y}}

\begin{example} \rm
\label{ex:counter}
Let $M = \R^2$ and $\lambda\colon ((a, b), (x, y)) \mapsto (x, y+a+bx)$. The orbits of
$\lambda$ are the lines parallel to the $y$-axis, so the foliate vector fields all have the form
$f(x)\dx + g(x,y)\dy$. 
The tangent vector fields have the form $g(x, y)\dy$.
However, the vector field $X = c(x, y) \dx +d(x, y) \dy$ is
 $G$-invariant iff $[\xi_M,X] = [\eta_M,X] = 0$, where $\xi_M=\dy$ and
$\eta_M=x\dy$ are the two generators of the action. This gives  
$$ \left[\xi_M,X\right] = \left[\dy,c\dx + d\dy\right] = c_y \dx + d_y \dy = 0,$$
implying $c_y = d_y = 0$, and
$$\left[\eta_M,X\right] = \left[x\dy,c\dx + d\dy\right] = -c \dx = 0,$$
implying $c = 0$. That is, $X = d(x) \dy$,
and the tangent and invariant vector
fields do not span the foliate vector fields.
\end{example}

The flow of a foliate vector field necessarily maps leaves to leaves diffeomorphically.  
In particular, each orbit is restricted to leaves of constant dimension.
However, for vector fields of the form $\Xt + \Xi$, even more is true.

\begin{theorem}
The flow of the G-foliate vector field $X = \Xt + \Xi$ preserves the isotropy
subgroup of the initial condition up to conjugacy. Specifically,
$G_{x(t)} = g(t) G_{x(0)} g(t)^{-1}$, where $g(t)$ satisfies Eq.
(\ref{eq:Gpiece}) below.
\end{theorem}
\begin{pf}
We write $x(t)=\lambda(g(t),h(t))$ as in Eq. (\ref{eq:representation}).
The flow $\phi_t$ of $X_{\rm inv}$
does not change the isotropy subgroup.  $G$-equivariance means
$\phi_t(\lambda(g,h(0)) = \lambda(g,\phi_t(h(0))$ for all $g\in G$,
so $g\in G_{h(0)}$ iff $g\in G_{h(t)}$. 
Then we have
$ G_{x(t)} = G_{g(t)h(t)} = g(t)G_{h(t)}g(t)^{-1} = g(t)G_{h(0)} g(t)^{-1}
= g(t)G_{x(0)} g(t)^{-1}$, 
as required.
\end{pf}

Note that the result is not true for all $G$-foliate flows. 
In Example \ref{ex:counter}, the
isotropy groups are $G_{(x,y)} = \{(a, b)\colon a = -bx\}$, which are all invariant under
conjugacy, because $G$ is abelian. Therefore $G_{(x_1,y)}$ is not conjugate to $G_{(x_2,y)}$
for $x_1 \ne x_2$. The
flow of the foliate vector field $\dot x  = 1$, $\dot y = 0$ therefore maps
leaves to leaves of non-conjugate isotropy.

\section{ Integrators for Lie group foliate vector fields}

Given a Lie group foliate vector field $X=X_{\rm tan} + X_{\rm inv}$,
this splitting into tangent and invariant pieces may not be unique, 
because it depends on the choice of $G$-invariant metric (or splitting
of $TM$).
In fact we have the whole family of splittings
$X = (X_{\rm tan}+Y) + (X_{\rm inv}-Y)$
for any $G$-invariant vector field $Y\in \Xtan$.
Nevertheless, once such a splitting has been found, it can be used
to construct foliate integrators. First, one can integrate each piece
separately and compose the results: $X_{\rm tan}$ by any integrator
for vector fields on homogeneous spaces \cite{mu-za}, and $X_{\rm inv}$
by any symmetry-preserving integrator, although these are 
easy to construct only when the action is linear. 

It is also possible to construct integrators in one piece.
Write the solution 
\begin{equation}
\label{eq:representation}
x(t)=\lambda(g(t),m(t)), \quad m(0)=x(0),\quad g(0)=1,
\end{equation}
and differentiate to get
$$ \dot x = D_1\lambda(g,m)\dot g + D_2\lambda(g,m)\dot m.$$
Here $D_1\lambda(g,m)\colon T_g G \to T_{x}M$ is the derivative
of $\lambda$ in its first slot, and $D_2\lambda(g,m)\colon T_{m}M \to T_{x}M$
is the derivative of $\lambda$ in its second slot.
Note that the first term on the right is tangent to the foliation.
Therefore, we choose 
$\dot g$ and $\dot m$ so that the first term is $X_{\rm tan}$ and
the second term is $X_{\rm inv}$:
\begin{eqnarray}
&&\dot g = (D_1\lambda(g,m))^{-1}X_{\rm tan}(\lambda(g,m)) \label{eq:Gpiece} \\
&&\dot m =(D_2\lambda(g,m))^{-1}X_{\rm inv}(\lambda(g,m)) = X_{\rm inv}(m).
\label{eq:Mpiece}
\end{eqnarray}
Here we are using that the range of $D_1\lambda(g,m)$ is $\Xtan(x)$.
If the action is not free then this map is not injective---many 
Lie algebra elements generate the same tangent vector $X_{\rm tan}(x)$---%
and some choice of the inverse $(D_1\lambda(g,m))^{-1}$ must be made to specify a
vector field on $G$. However, in the examples below a natural choice can be
made. Thus we have extended the vector field on $M$ to a vector field
on $G\times M$ with extra foliate structure of its own. Integrating
the 
reduced vector field $X_{\rm inv}(m)$ by a symmetry-preserving
method leaves a reconstruction problem on $G$. However, the vector
field on $G$ is now known only at the time steps of the solution of 
$X_{\rm inv}(m)$ which could make it difficult to achieve high orders.
It is better to integrate the full system together as follows.

\begin{theorem} When $M=\R^n$ and the group action is linear, 
the following algorithm provides
a foliate integrator of order $p$ for $X=X_{\rm tan} + X_{\rm inv}$.
(i) Let $x(t) = \lambda(g(t),m(t))$, $g(0)=1$, $m(0)=x(0)$ and
choose a differential equation on $G$ according to Eq. (\ref{eq:Gpiece}).
(ii) Apply a Runge-Kutta method of order $p$ 
to $\dot m = X_{\rm inv}(m)$, and the
associated Runge-Kutta--Munthe-Kaas (RKMK) \cite{mu-za} method 
to the chosen differential equation on $G$, with stage values for
$m$ given by the corresponding stage values of the RK method applied
to $\dot m = X_{\rm inv}(m)$.
\end{theorem}
\begin{pf}
Since the group action is linear, $X_{\rm inv}$ has a linear symmetry
and this is preserved by the Runge-Kutta method. That is, the integrator
on $m$ is $G$-equivariant and hence preserves the foliation into $G$-orbits.
Applying RKMK to the differential equation on $G$ is equivalent to
applying the same RKMK method to the full system on $G\times M$
with respect to the group action of $G\times M$ on itself given by
$$\lambda((g_1,m_1),(g_2,m_2)) = (g_1g_2,m_1+m_2)$$
and hence it has order $p$. Reconstructing the solution by acting on
$m(t)$ with $g(t)$ only moves the point around on the new leaf, so the
total integrator is foliate.
\end{pf}

\begin{example} \rm
The matrix differential equation $\dot L = [A(L),L]$ motivated much of
the work in Lie group integrators \cite{iserles,ca-is-za}. Let $
G\subset GL(n)$ be a matrix
Lie group, let $M=\mathfrak{g}$, the Lie algebra of $G$, and let $G_1$
be a subgroup of $G$ which acts on $M$ by adjoint action, i.e. 
$\lambda(U,L) = U L U^{-1}$ for $U\in G_1$, $L\in \mathfrak{g}$.
Recall that the ``isospectral manifolds'' of $\mathfrak{gl}(n)$ are the sets of
matrices similar by an element of $GL(n)$. In the present
example, the
leaves of the foliation 
are the sets of matrices in $\mathfrak{g}$
which are similar by an element of $G_1$, and hence are submanifolds
of the isospectral manifolds. 

All vector fields tangent to the leaves can be written in 
the form $[A(L),L]$, where $A\colon \mathfrak{g}\to\mathfrak{g}_1$. From Theorem
\ref{theorem:main},
when $\mathfrak{g}$ admits an adjoint invariant metric (e.g. when
it is compact), 
all foliate vector fields can be written in the form
\begin{equation}
\label{eq:15}
\dot L = [A,L] + f(L),\quad f(ULU^{-1}) = U f(L) U^{-1}\ \forall U\in G_1
\end{equation}
for some $(G_1)$-adjoint-invariant function $f\colon \mathfrak{g}\to\mathfrak{g}$.

It is not so easy to explicitly construct all invariant vector fields. The classical
method requires determining a complete and independent set of differential
invariants of the action. Another approach which can work if the action is
understood well enough is to choose a reference point on each orbit, choose
an isotropy-invariant tangent vector at each reference point, and push them     
around the leaves using the group action.

For example, consider the case $G = G_1 = GL(n)$. The orbits are the conjugacy     
classes of matrices so a reference point on each orbit is given by the 
Jordan matrices.
Restricting for simplicity to the diagonalizable matrices, whose isotropy   
groups are the group of nonsingular diagonal matrices, we first choose $f(\Lambda)$  
diagonal (so that it is isotropy invariant), and then define
$$f(L) = Ug(\Lambda)U^{-1}$$
where $L = U\Lambda U^{-1}$ is the diagonalization of $L$. This constructs all invariant
vector fields on the diagonalizable matrices.

Furthermore, this case provides a nice example of Theorem 8 in a case in
which the action is not proper. The set of diagonal matrices form a slice at
any diagonal matrix, on which all isotropy subgroups are equal, so assumption
(ii) of the theorem is automatically satisfied. The $G_{x_0}$-invariant splitting of $\g$
(assumption (i)) is provided by splitting into diagonal and off-diagonal parts.
Hence from Theorem 8, Eq. (\ref{eq:15}) spans all foliate vector fields in this case. We 
can think of the foliate vector fields as those for which the evolution of the
eigenvalues of $L$ depends only on those eigenvalues, and not on the individual
entries of $L$.

More explicit examples of invariant vector fields are given by
$$ f(L) = p(L)g(\tr L,\tr L^2,\dots,\tr L^n),$$
where $p$ is a real analytic function (which, when extended to matrices, maps
$\g$ into $\g$)and $g\colon \R^n\to \R$.

To construct foliate integrators, following Eq. (\ref{eq:representation}) 
we represent the
solution as $L(t) =U(t) F(t) U(t)^{-1}$, $U(t)\in G_1$, $F(t)\in \mathfrak{g}$,
and work with the system on $G_1\times\mathfrak{g}$,
$$\eqalign{
\dot U &= A(UFU^{-1})U,\quad U(0)=I,\cr
\dot F &= f(F), \quad F(0)=L(0),\ \cr
}$$
Since the symmetry of the second equation is linear, it is preserved
by any linearly covariant method, such as Runge-Kutta. Runge-Kutta--Munthe-Kaas
applied to the whole system provides a foliate integrator.
The simplest example
is the following first-order ``Lie-Euler" method. Noting that we can take $U_n =
I$, $F_n = L_n$ at each time step (only calculating the updates which move $L$
around on the leaves), we get the method
$$U_{n+1} = \exp(\tau A(F_n)),\quad  F_{n+1} = F_n + \tau f(F_n),$$
or, in terms of $L$,
$$L_{n+1} = \exp(\tau A(L_n))(L_n + \tau f(L_n)) \exp(-\tau A(L_n)).$$

\end{example}

\begin{example} \rm
Let $M=\R^{n\times p}$ and let the matrix Lie group $G\subset GL(n)$
act on $M$ by left multiplication, $\lambda(U,A) = UA$ for
$U\in G$, $A\in M$. 
As in the previous
example, in the case that the action is proper (e.g. when $G$ is compact) we 
can write the foliate vector fields in the form
$$ \dot A = g(A)A + f(A)$$
for some functions $g\colon M\to\mathfrak{g}$ and $f\colon A\to A$, where $f$ is
invariant, i.e. $f(UA) = Uf(A)$ for all $U\in G$. 
As before it is difficult to find all invariant vector fields.
Examples are $f(A) = AV (A^{\rm T}A)$ for any $V\colon \R^{p\times p}
\to\R^{p\times p}$, for $G = SO(n)$,
and $f(A) = AV (\det A)$ for any $V\colon \R \to \R^{p\times p}$ for $G = SL(n)$. 
The planar
systems of Eqs. (\ref{eq:polar1},\ref{eq:polar2}) are examples with $n = 2$, $p = 1$, and $G = SO(2)$. For $G =       
SO(n)$, the leaves are the sets $A^{\rm T}A = $ const., one of which is the Stiefel
manifold $A^{\rm T}A = I$.

To construct foliate integrators, following Eq. (\ref{eq:representation}) 
we represent
the solution as $A(t) = U(t)F(t)$, $U(t)\in G$, $F(t)\in R^{n\times p}$,
and work with the system 
$$\eqalign{
\dot U &= g(UF)U,\quad U(0)=I, \cr
\dot F &= f(F), \quad F(0)=A(0),\ \cr
}$$
As in the previous example, RKMK applied to the whole system
provides a foliate integrator. The Lie--Euler method is
\begin{equation}
\label{eq:lie-euler}
A_{n+1} = \exp(\tau g(A_n))(A_n + \tau f(A_n))\end{equation}
An example of this method is shown in Figure 2.

\end{example}

\begin{figure} \begin{center}
\vbox{\includegraphics[width=8cm,keepaspectratio,clip]{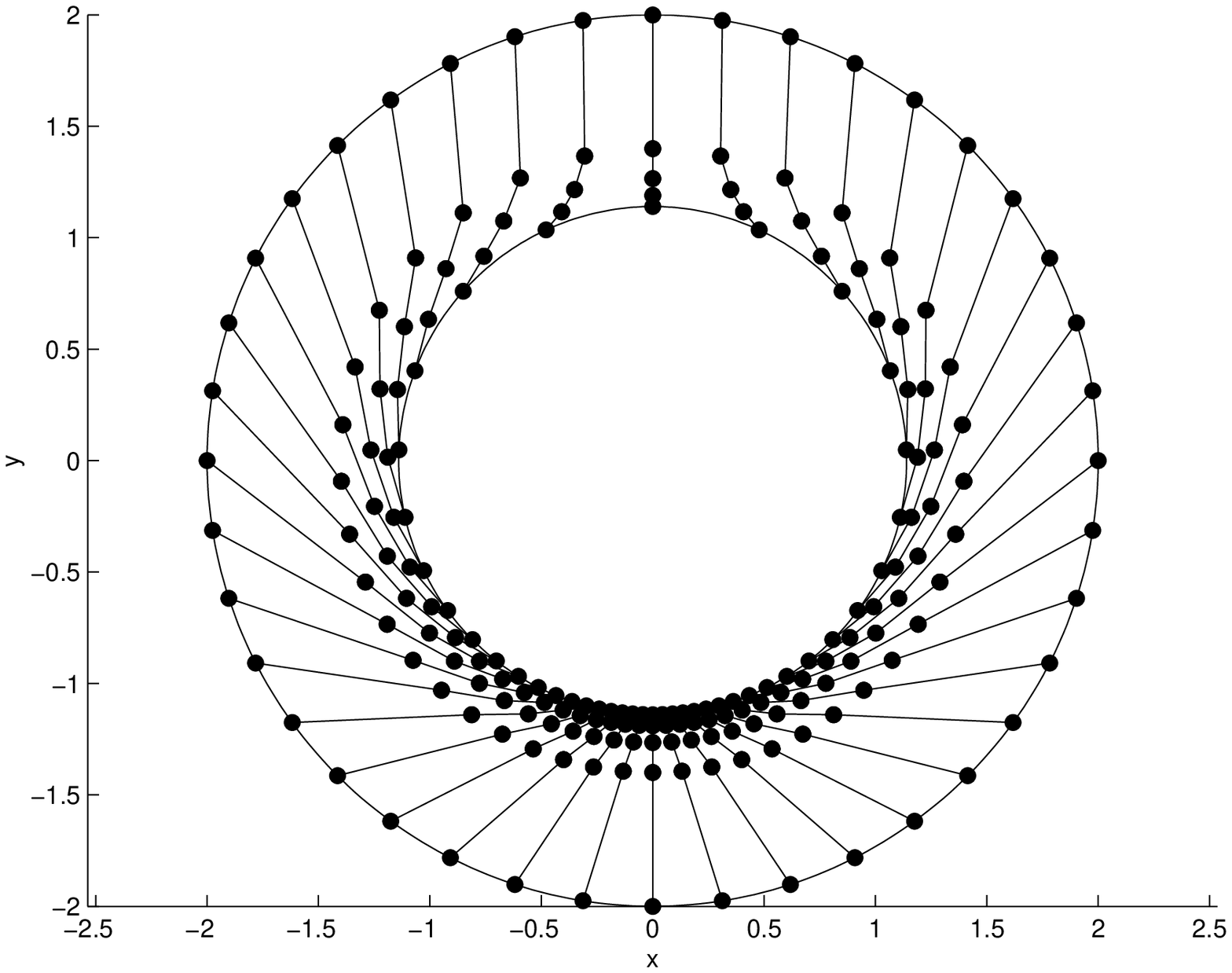} }
\bigskip
\vbox{\includegraphics[width=8cm,keepaspectratio,clip]{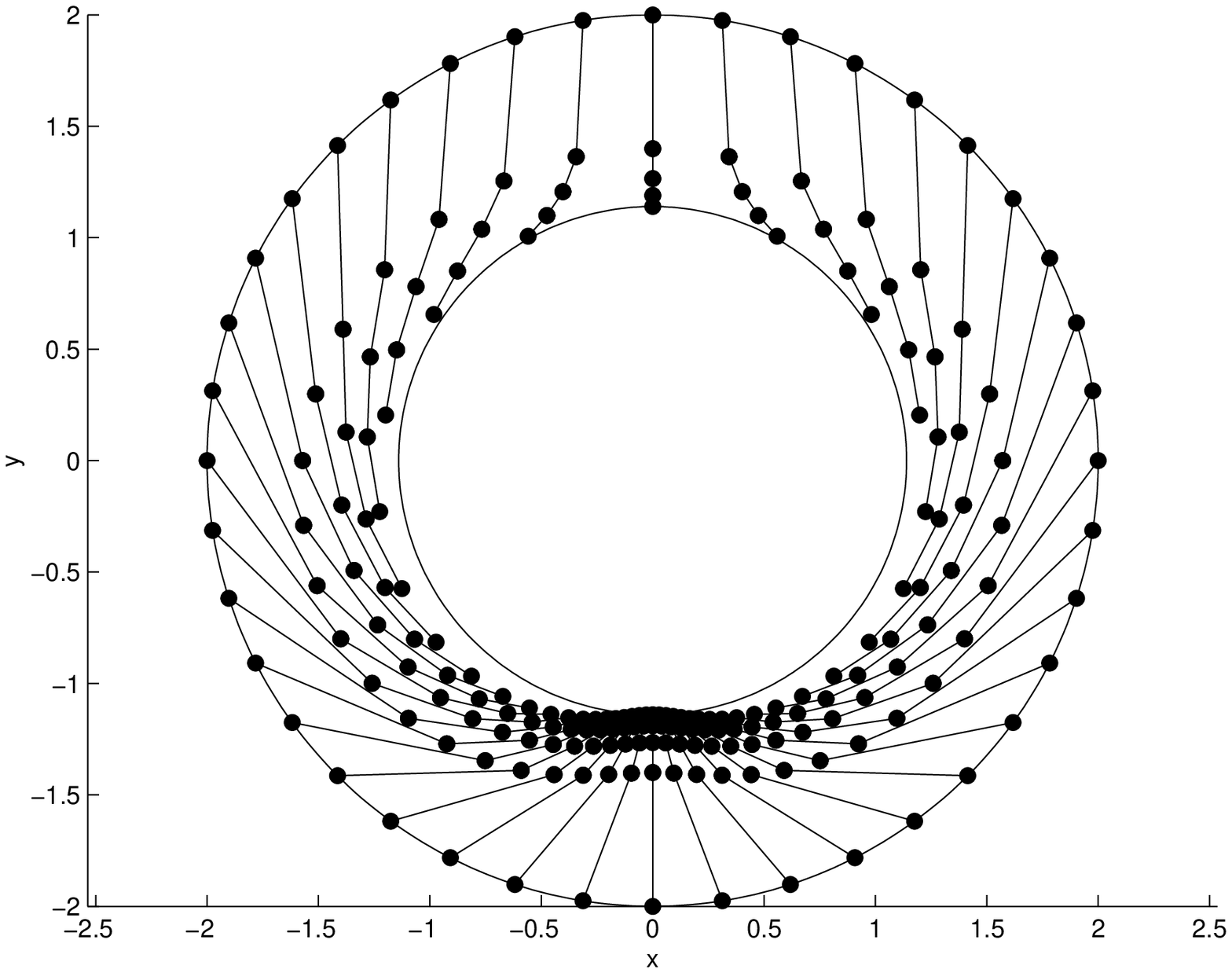} }
\caption{Foliate vs. nonfoliate integrators. The ODE of Eq. (1) is integrated by the
(foliate) Lie-Euler method Eq. (\ref{eq:lie-euler}) (top), 
and the (nonfoliate) Euler method (bottom).
The 20 initial conditions lie on a circle of radius 2, and four time steps of 0.1
are shown. In the nonfoliate integrator, the final values do not lie on the reference
circle shown. Note that the two methods coincide on $x = 0$, where the tangential
component vanishes.}
\end{center}
\end{figure}

\section{Conclusions}

The situation considered here, of 
vector fields preserving a given foliation, can be extended in
several ways. First, a system may preserve several different foliations.
This is formalized in the `multifoliate' structure introduced by Kodaira
and Spencer \cite{ko-sp}: the distributions form a lattice, closed under intersections
and (an appropriate) join. A simple example on $\R^3$ is provided by the two
distributions $\dx$ and $\dy$,
for which the multifoliate vector fields have the form
$$\dot x = f(x,z),\quad \dot y =g(y,z),\quad \dot z = h(z).$$
A more sophisticated case occurs in Hamiltonian systems with symmetry,
which preserve the foliations into group orbits and into momentum level sets,
and in Poisson and conformal Poisson systems, which preserve in addition the
foliation into symplectic leaves \cite{mc-pe}. We plan to study such systems more in
the future \cite{mc-pe2}.

Second, a foliate system may have extra structure, corresponding to
the different infinite dimensional Lie subalgebras of the Lie algebra
of foliate vector fields. These so-called nonprimitive Lie subalgebras
have not been classified. However, many examples can be constructed by
considering (i) the vector field to lie in some other Lie algebra, as
of Hamiltonian or volume-preserving vector fields; (ii) the reduced
system on the space of leaves to lie in some other Lie algebra (in this
case the foliated phase space has {\em transverse} structure); and
(iii) the reconstruction system, considered as a nonautonomous vector
field on a leaf, to lie in some other Lie algebra. (It is the relationships
between the reduced and reconstruction systems 
which have not been classified.) Cartan
[5] provides an interesting list of all infinite-dimensional Lie pseudogroups in
two dimensions. For example, in the group
$$ (x, y) \mapsto \left(\frac{  xf(y) + \phi(y)}{x\psi(y)+1},g(y)\right),$$
the foliation $y =$ const. is preserved, and the reconstruction dynamics lie in
the fractional linear group. But in the group
$$(x,y)\mapsto \left(x f'(y)^{-2} - f'''(y)f'(y)^{-2} + 
\frac{3}{2}f''(y)^2 f'(y)^{-4},f(y)\right),$$
what structure is preserved? In all of these cases it makes sense to ask how the
structure affects the dynamics (e.g., what are its homeomorphism invariants),
how it can be detected in a given system, and how integrators that lie in the   
corresponding group of diffeomorphisms can be constructed. From this point   
of view the results presented here are just a beginning.

\begin{ack}
We would like to thank the Marsden Fund of the Royal Society of New Zealand,    
the Australian Research Council, and the Center for Advanced Studies, Oslo,   
for their support.
\end{ack}

\end{document}